# ON GENERALIZATIONS OF THE BUTTERFLY THEOREM

NICHOLAS PHAT NGUYEN

**Abstract**. This paper provides a uniform explanation of different extensions and generalizations of the butterfly theorem based on the Desargues' involution theorem.

1.  **INTRODUCTION**. The butterfly theorem is a perennial fascination in Euclidean plane geometry. There have been many extensions and generalizations of the butterfly theorem, some of them published in recent years. In this paper, we will show how these extensions and generalizations are just different manifestations of the Desargues' involution theorem.

While it is known that the original butterfly theorem is a special case of the Desargues' involution theorem, it is less clear what the connection is between the Desargues' involution theorem and different extensions and generalizations of the butterfly theorem that we have seen over the years. The published proofs for these generalizations use a variety of methods, all rather complicated. The paper will show how these different generalizations of the butterfly theorem follow in a natural and rather simple way from the Desargues' involution theorem.

For the convenience of the readers, we recall the Desargues involution theorem below.

———————————————————————————.

**Theorem (Desargues' Involution Theorem) – Real projective plane version:** *Consider four points in general position in the real projective plane, i.e., no three of these four points are collinear. Let $\mathcal{F}$ be the family of conics passing through these four base points. Then for any line $\ell$ that does not pass through any of the four base points, each conic in the family $\mathcal{F}$ will, if it intersects with the line $\ell$, do so in a pair of points that are conjugate under an involution of the line $\ell$.*

Here, an involution means a projective transformation of the line $\ell$ that has order two. In order words, it refers to a projective transformation $T$ of $\ell$ such that $T \neq Id$, and $T^2 = Id$. The theorem says that there is an involution $T$ such that for any conic $G$ in the family $\mathcal{F}$, we have one of the following situations: (i) the conic $G$ does not intersect $\ell$ at all; or (ii) the conic $G$ intersects $\ell$ in two distinct points $x$ and $y$ where $y = T(x)$; or (iii) the conic $G$ intersects $\ell$ in a single point $z$, where $z = T(z)$, i.e., $z$ is a fixed point of the involution $T$.

An involution of a projective line either has no fixed point or exactly two fixed points. If $I$ and $J$ are two fixed points of an involution, then any pair $(U, V)$ of conjugate points of the involution have the property that the cross ratio $(I, J; U, V) = -1$, i.e., the points $I, J, U, V$ form a harmonic range. Hence an involution with two given fixed points is uniquely determined.

Note that the Desargues' involution theorem also holds over the complex numbers. Because a conic will always intersect a given line in a complex projective plane, and because an involution of a complex projective line always has two fixed points, we have the following version of the Desargues' involution theorem over the complex numbers.

**Theorem (Desargues' Involution Theorem) – Complex projective plane version:** *Consider four points in general position in the complex projective plane. Let $\mathcal{F}$ be the family of conics passing through these four points. Then for any line $\ell$ that does not pass through any of these four base points, each conic in $\mathcal{F}$ will intersect the line $\ell$ either in a pair of distinct points that are conjugate under an involution of the line $\ell$, or in a single point that is one of the two fixed points of the involution.*

**2. GENERALIZATIONS OF THE BUTTERFLY THEOREM.** We use the same set-up and notation as in the Desargues' involution theorem stated above. We will refer to the involution on the line $\ell$ induced by the family $\mathcal{F}$ as the Desargues involution.

**Proposition 1:** *Suppose a conic in the family $\mathcal{F}$ intersects the line $\ell$ in a single point M, and another conic in the family $\mathcal{F}$ intersects the line $\ell$ in two distinct points P and Q. Let N be the unique point of the line $\ell$ such that the cross ratio $(M, N; P, Q) = -1$. Then for a conic in the family $\mathcal{F}$ we have one of the following situations: (i) the conic does not intersect the line $\ell$; or (ii) the conic intersects the line $\ell$ in a single point which is either M or N, or (iii) the conic intersects the line $\ell$ in a pair of distinct points S and T such that $(M, N; S, T) = -1$.*

*Proof.* Consider the unique involution of the line $\ell$ with $M$ and $N$ as the fixed points. This involution agrees with the action of the Desargues involution on the points $M, P$ and $Q$. Since a projective transformation of a projective line is completely determined by its effect on three distinct points, these two involutions must be the same. ∎

In particular, if we choose a coordinate on the projective line $\ell$ so that the point $M$ is the origin (coordinate zero) and $N$ is the point at infinity, then the coordinates of any pair of conjugate points are $x$ and $-x$, i.e., $M$ is the mid-point of each pair of distinct conjugate points.

The classic butterfly theorem is an example of the above basic butterfly configuration, where all four base points $A, B, C, D$ lie on a circle that intersects the line $\ell$ in two distinct points $P$ and $Q$, and the chords $AC$ and $BD$ intersect the line $\ell$ at the midpoint $M$ of $P$ and $Q$. The generalized butterfly theorem has the same basic configuration, but uses a proper conic instead of a circle. See [2], [3] and [4].

We could also have a line $\ell$ passing through a diagonal point $M$ of the quadrangle $ABCD$ and two opposite sides of that quadrangle intersecting the line $\ell$ in two distinct points symmetric about M. Another variant configuration is to use a tangent line to a conic rather than a chord.

**Proposition 2:** *Suppose that two different conics in the family $\mathcal{F}$ intersect the line $\ell$ in two different pairs of (real or imaginary) distinct points (P, Q), and (S, T) over the complex numbers. Also, suppose there are two distinct points M and N on the line $\ell$ with cross ratio (M, N; P, Q) = −1 = (M, N; S, T). In that case, we have one of the following situations for any conic in the family $\mathcal{F}$: (i) the conic does not intersect the line $\ell$, in which case it will intersect $\ell$ in two imaginary points U and V (over the complex numbers) with cross ratio (M, N; U, V) = −1; or (ii) the conic intersects the line $\ell$ at a single point, which is either M or N, or (iii) the conic intersects the line $\ell$ at two distinct points U and V such that the cross ratio (M, N; U, V) = −1.*

*Proof.* If we extend the ground field to the complex numbers, then we have a Desargues involution on the complex version $\ell(\mathbb{C})$ of $\ell$. Consider the unique involution of the complex projective line $\ell(\mathbb{C})$ with M and N as the fixed points. This involution agrees with the action of the Desargues involution on the points P, Q, S and T. Hence it must be the Desargues involution on $\ell(\mathbb{C})$.

If we restrict the Desargues involution of $\ell(\mathbb{C})$ to the real points of that line, we get the Desargues involution on the real projective line $\ell$ that follows from the real version of the Desargues involution theorem. This means the Desargues involution on the real projective line $\ell$ must be the unique involution with M and N as the fixed points. ∎

If we choose a coordinate on the projective line $\ell$ so that the point M is the origin (coordinate zero) and N is the point at infinity, then M is the mid-point of each pair of distinct conjugate points.

Many generalizations of the butterfly theorem are examples of Proposition 2. The generalization by Murray Klamkin in [2] is the case where we have two conics in the family $\mathcal{F}$ intersecting a given line $\ell$ in two pairs of distinct points with the same mid-point M.

As another example, consider the case where the four base points A, B, C, and D lie on a circle in a Euclidean plane. Suppose that a conic passing through these four points happens to intersect a given line $\ell$ at two points $P$ and $Q$ that are equidistant from the center of the circle. Let $M$ be the mid-point of $P$ and $Q$. The condition that $P$ and $Q$ are equidistant from the center of the circle (which is equivalent to the line from the center to the mid-point $M$ being orthogonal to the line $\ell$) means that if we allow complex numbers then the circle will intersect the line $\ell$ in two (real or imaginary) points whose mid-point is also $M$. It follows from Proposition 2 that any conic passing through the base points $A$, $B$, $C$, and $D$ will either intersect the line $\ell$ at a single point which is $M$ or the point at infinity, or it will intersect the line $\ell$ in two distinct (real or imaginary) points whose mid-point is $M$. This is the generalization of the butterfly theorem stated in [5], [7], and [11].

Incidentally, in this configuration of four base points on a circle and a conic passing through the base points and intersecting a given line $\ell$ at two points $P$ and $Q$ that are equidistant from the center of the circle, Proposition 2 shows that the mid-point $M$ of $P$ and $Q$ and the point $N$ at infinity are the fixed points of the Desargues' involution on the line $\ell$. That means the family $\mathcal{F}$ of all conics passing through the 4 base points contains: (i) a conic that is tangent to the line $\ell$ at $M$, and (ii) a conic that is tangent to the line $\ell$ at $N$, i.e., a hyperbola with the line $\ell$ as one asymptote. This is a quite interesting fact by itself, which does not seem to be noted in the literature.

For a more general example, consider the case where: (a) a conic $G$ in the family $\mathcal{F}$ intersects the line $\ell$ in two distinct points $P$ and $Q$; and (b) the line $\ell$ happens to be conjugate to a diameter $k$ of another conic $H$ in the family $\mathcal{F}$ and the lines $\ell$ and $k$ intersect at the mid-point $M$ of $P$ and $Q$. Here, being a diameter of $H$ means that the line $k$ passes through the pole (relative to the conic $H$) of the line at infinity. That pole is also known as the center of the conic $H$.

For the line $\ell$ to be conjugate to $k$ means that the line $\ell$ passes through the pole $N$ of the line $k$, and that pole is on the line at infinity because the line $k$ goes through the center of the conic $H$. In this situation, the line $\ell$ will intersect the conic $H$ (over the complex numbers) in two distinct points $S$ and $T$ such that $(M, N, S, T) = -1$. If we chose a coordinate on the line $\ell$ so that $M$ is at the origin and $N$ is the point at infinity, then $M$ is the mid-point of $P$ and *Q* as well as *S* and *T*. Proposition 2 then applies in this case and shows us that the conics in the family $\mathcal{F}$ will generally intersect the line $\ell$ (when they do) in pairs of points symmetric about *M*. This is the generalization of the butterfly theorem contained in [9].

An example of this situation is when the line $\ell$ happens to be perpendicular to an axis of a conic in the family $\mathcal{F}$. An axis is one of two diameters that passes through 2 opposite foci of the conic. If $\ell$ happens to be perpendicular to such an axis, then it passes through the pole of that axis on the line at infinity and hence is conjugate to the axis. This is the generalization of the butterfly theorem contained in [8].

3.  **BUTTERFLY POINTS.**   In general, if we can find a point *M* in the affine plane distinct from the four base points and a point *N* on the line at infinity that are conjugate to each other relative to two different conics in the family $\mathcal{F}$, then we have a butterfly configuration. Indeed, in that case the points *M* and *N* will be conjugate to each other relative to all the conics in the family $\mathcal{F}$. That means the conics in the family $\mathcal{F}$ will generally intersect the line passing through *M* and *N* (when they do) in pairs of points that form a harmonic range with *M* and *N*, i.e., points that are symmetric about *M*. The exceptions are the two conics in the family $\mathcal{F}$ that are tangent to the line at either *M* or *N*.

How do we find such a point *M*, which we will refer to as a butterfly point of the four given base points? Given an arbitrary point *M* distinct from the four base points, it is well-known that all the polars of *M* relative to the conics in the family $\mathcal{F}$ constitute a pencil of lines passing through a common point. See [1] at p. 342. So if there is a point *N* on the line of infinity that is conjugate to *M* relative to all the conics in the family $\mathcal{F}$, then we know *N* must be that common point. Moreover, the line at infinity must also be the polar of *M* relative to some conic in the family $\mathcal{F}$. So *M* must be the center of a conic in $\mathcal{F}$.

Conversely, if *M* is the center of some conic in the family $\mathcal{F}$, then the line at infinity is by definition the polar of *M* relative to that conic. The polar of *M* relative to another conic in the family $\mathcal{F}$ will be a different line which will intersect the line at infinity in a point *N*. The points *M* and *N* are then conjugate to each other relative to two different conics in $\mathcal{F}$, and we have a butterfly configuration with *M* as the butterfly point.

Therefore all the butterfly points relative to the four base points are the centers of all the conics in the family $\mathcal{F}$. This is stated in [10] as its generalization of the butterfly theorem. These centers, the poles of the line at infinity relative to all the conics in $\mathcal{F}$, is a well-known conic known as the eleven-point conic defined by the four given base points. This eleven-point conic passes through 11 points in the complex projective plane defined by the configuration of the four base points A, B, C, D: the 3 diagonal points of the quadrangle ABCD, the 6 mid-points of the sides of the quadrangle ABCD, and the 2 points on the line at infinity that are the fixed points of the Desargues involution induced by the family $\mathcal{F}$. See [1] at p. 342-343. If the four base points lie on a circle, then this eleven-point conic also passes through the center of the circle (because it is the locus of all centers of conics in the family $\mathcal{F}$). Moreover, because the two fixed points of the Desargues involution on the line at infinity form a harmonic range with the two intersection points of the circle and the line at infinity (the imaginary circular points at infinity), the eleven-point conic in this case is actually a rectangular hyperbola.

NICHOLAS PHAT NGUYEN
*12015 12th Dr. SE, Everett, WA 98208, U.S.A*
*Email: nicholas.pn@gmail.com*